# Random walk in random scenery:
# A survey of some recent results

## Frank den Hollander[1,2,*] and Jeffrey E. Steif[3,†]

*Leiden University & EURANDOM and Chalmers University of Technology*

**Abstract.** In this paper we give a survey of some recent results for random walk in random scenery (RWRS). On $\mathbb{Z}^d$, $d \geq 1$, we are given a random walk with i.i.d. increments and a random scenery with i.i.d. components. The walk and the scenery are assumed to be independent. RWRS is the random process where time is indexed by $\mathbb{Z}$, and at each unit of time both the step taken by the walk and the scenery value at the site that is visited are registered. We collect various results that classify the ergodic behavior of RWRS in terms of the characteristics of the underlying random walk (and discuss extensions to stationary walk increments and stationary scenery components as well). We describe a number of results for scenery reconstruction and close by listing some open questions.

## 1. Introduction

Random walk in random scenery is a family of stationary random processes exhibiting amazingly rich behavior. We will survey some of the results that have been obtained in recent years and list some open questions. Mike Keane has made fundamental contributions to this topic. As close colleagues it has been a great pleasure to work with him. We begin by defining the object of our study.

Fix an integer $d \geq 1$. Let $X = (X_n)_{n \in \mathbb{Z}}$ be a sequence of i.i.d. random variables taking values in a *possibly infinite* set $F \subset \mathbb{Z}^d$ according to a common distribution $m_F$ having full support on $F$. Let $S = (S_n)_{n \in \mathbb{Z}}$ be the corresponding two-sided RANDOM WALK on $\mathbb{Z}^d$, defined by

$$S_0 = 0 \qquad \text{and} \qquad S_n - S_{n-1} = X_n, \, n \in \mathbb{Z},$$

i.e., $X_n$ is the step at time $n$ and $S_n$ is the position at time $n$. To make $S$ into an *irreducible* random walk, we will assume that $F$ generates $\mathbb{Z}^d$, i.e., for all $x \in \mathbb{Z}^d$ there exist $n \in \mathbb{N}$ and $x_1, \ldots, x_n \in F$ such that $x_1 + \cdots + x_n = x$. The *simple random walk* is the case where $F = \{e \in \mathbb{Z}^d : |e| = 1\}$ and $m_F(e) = \frac{1}{2d}$ for $e \in F$.

*The research of FdH was supported by the Deutsche Forschungsgemeinschaft and the Netherlands Organisation for Scientific Research through the Dutch-German Bilateral Research Group on "Mathematics of Random Spatial Models from Physics and Biology".

† The research of JES was supported by the Swedish National Science Research Council and by the Göran Gustafsson Foundation (KVA).

[1]Mathematical Institute, Leiden University, P.O. Box 9512, 2330 RA Leiden, The Netherlands, e-mail: denholla@math.leidenuniv.nl

[2]EURANDOM, P.O.Box 513, 5600 MB Eindhoven, The Netherlands.

[3]Department of Mathematics, Chalmers University of Technology, Gothenburg, Sweden, e-mail: steif@math.chalmers.se

Next, let $C = (C_z)_{z \in \mathbb{Z}^d}$ be i.i.d. random variables taking values in a *finite* set $G$ according to a common distribution $m_G$ on $G$ with full support. Unless stated otherwise, we will restrict to the case where $G = \{-1, +1\}$ and $m_G(-1) = m_G(+1) = \frac{1}{2}$, although most results in this paper hold in general. We will refer to $C$ as the RANDOM SCENERY, i.e., $C_x$ is the scenery value at site $x$.

In what follows, $X$ and $C$ will be taken to be *independent*. Let

$$Y = (Y_n)_{n \in \mathbb{Z}} \qquad \text{with} \qquad Y_n = (C \circ S)_n = C_{S_n}$$

be the sequence of scenery values observed along the random walk. We will refer to $Y$ as the SCENERY RECORD. The joint process

$$Z = (Z_n)_{n \in \mathbb{Z}} \qquad \text{with} \qquad Z_n = (X_n, Y_n)$$

is called RANDOM WALK IN RANDOM SCENERY (RWRS). This process registers both the step taken by the walk and the scenery value at the site that is visited. Note that, while $X$ is simple, $Y$ is complicated (because it is a composition of two random processes). The interplay between $X$ and $Y$ will be important.

We will assume that the reader is familiar with a number of key concepts from ergodic theory, namely, K-AUTOMORPHISM, BERNOULLI, WEAK BERNOULLI, and FINITARY FACTOR OF AN I.I.D. PROCESS. For definitions we refer to Walters [41]. For reasons of exposition, we give loose definitions here. A K-automorphism is a random process with a trivial future tail $\sigma$–field. A Bernoulli process is one that can be coded (in an invertible and time-invariant manner) from an i.i.d. process. A weak Bernoulli process is one where the past and the far distant future have a joint distribution that is close in total variation norm to the distribution of the past and the far distant future put together independently. A finitary factor of an i.i.d. process is a Bernoulli process for which the coding is such that, in order to determine one of the output bits, one need only look at a finite (but random) number of bits in the i.i.d. process.

In ergodic theory, when $d = 1$ and the random walk is simple, RWRS is referred to as the $T, T^{-1}$-process. The reason is that if $T$ is the shift on the scenery sequence, then with each step the walker sees the scenery sequence shifted either by $T$ or by $T^{-1}$ depending on whether the step is to the right or to the left. The interest in RWRS originally came from the fact that, for simple random walk in $d = 1$, $Z$ was conjectured to be a natural example of a K-automorphism that is not Bernoulli. As the history given below reveals, this conjecture turned out to be true. In our opinion it is *by far the simplest* such example (in terms of the description of the process, though not in terms of the proof).

The outline of this paper is as follows. In Section 2 we focus on the ergodic properties of RWRS and present a history of results that have been obtained so far, organized in Sections 2.1-2.6. In Section 3 we describe a number of results for scenery reconstruction. In Section 4 we close by listing some open questions. Inevitably, what follows is *our* selection of highlights and omits certain works that could also have been included.

## 2. Ergodic properties

In Sections 2.1 and 2.2 we list the main theorems that determine when RWRS is a K-automorphism, is Bernoulli, is weak Bernoulli, or is a finitary factor of an i.i.d. process. In Section 2.3 we have a look at when these properties hold for reduced RWRS, the second component of RWRS, i.e., the scenery record alone. In Section 2.4 we make a brief excursion into random walks with stationary increments and random sceneries with stationary components, in order to see what properties survive when we relax the i.i.d. assumptions. In Section 2.5 we consider induced RWRS, which is obtained by observing the scenery only when a +1 is visited. Finally, in Section 2.6 we investigate the continuity properties of the conditional probability distribution for RWRS at time zero given the configuration at all other times.

### 2.1. K-automorphism, Bernoulli and weak Bernoulli

RWRS is clearly stationary. Since the walk may return to sites visited before, it is not i.i.d. By Kakutani's random ergodic theorem, it is ergodic. Our starting point is the following general result.

**Theorem 2.1 (Meilijson [33]).** *RWRS associated with an arbitrary random walk is a K-automorphism (i.e., has a trivial future tail $\sigma$–field).*

The intuition behind this result is that the distribution of the walk spreads out for large times, so that the walk sees an ergodic average of the scenery. Although Theorem 2.1 was proved only for $d = 1$ (in the more general setting of so-called skew-products), the argument easily extends to arbitrary $d \geq 1$. See Rudolph [36] for related results.

We point out that for the case of simple random walk in $d = 1$, the result in Theorem 2.1 was "known" prior to [33]. On the other hand, Meilijson actually proved his result for any totally ergodic random scenery (See Section 2.4 below).

It was known early on that being Bernoulli implies being a K-automorphism, because the latter is isomorphism invariant (i.e., invariant under coding) and an i.i.d. random process has a trivial future tail $\sigma$-field. In 1971, Adler Ornstein and Weiss conjectured that RWRS associated with simple random walk in $d = 1$ is not Bernoulli (see [42] and [16]). If true, then this would provide a beautiful and natural example of a K-automorphism that is not Bernoulli. (At some earlier stage, it was an open question whether every K-automorphism was Bernoulli. Counterexamples were constructed by Ornstein and later by Ornstein and Shields, but these were much less natural.) The conjecture was settled in a deep paper by Kalikow.

**Theorem 2.2 (Kalikow [16]).** *RWRS associated with simple random walk in $d = 1$ is not Bernoulli.*

(In fact, Kalikow actually proves the stronger result that the process is not "loosely Bernoulli", a notion we will not consider.)

Theorem 2.2 was later extended to cover an almost arbitrary recurrent random walk.

**Theorem 2.3 (den Hollander and Steif [13]).** *If the random walk is recurrent with $\sum_{x \in \mathbb{Z}^d} |x|^\delta m_F(x) < \infty$ for some $\delta > 0$ and satisfies a certain technical condition, then the associated RWRS is not Bernoulli.*

We will not explain the "technical condition", because it is extremely weak and is in fact conjectured in [13] to hold for an arbitrary random walk (!). In [13] it is proved that if $m_F$ has one component that is in the domain of attraction of a stable law with index $\geq 1$, then the technical condition is already fulfilled.

Theorem 2.3 shows that recurrence of the random walk essentially implies that the associated RWRS is not Bernoulli. The next result tells us that if the random walk is transient, then the associated RWRS is Bernoulli, providing us with a nice dichotomy.

**Theorem 2.4 (den Hollander and Steif [13]).** *If the random walk is transient, then the associated RWRS is Bernoulli.*

The concept of weak Bernoulli was studied in the early days of ergodic theory. It was introduced by Kolmogorov under the name absolutely regular, and is called $\beta$–mixing in some circles of probabilists. Weak Bernoulli was known early on to imply Bernoulli. The fact that the two are not equivalent was proved by Smorodinsky [38]. The characterization of when RWRS is or is not weak Bernoulli turns out to be very interesting. Rather than describing the results in complete generality, we restrict to special classes of random walk.

**Theorem 2.5 (den Hollander and Steif [13]).**
(i) *For simple random walk on $\mathbb{Z}^d$, the associated RWRS is weak Bernoulli if and only if $d \geq 5$.*
(ii) *For random walk with bounded step size and nonzero mean, the associated RWRS is weak Bernoulli for any $d \geq 1$.*

The proof of Theorem 2.5 is based on a coupling argument where, given two independent pasts of RWRS, the futures are coupled so as to make them agree far out. Loosely speaking, weak Bernoulli is equivalent to such a coupling being possible.

The *phase transition* in the behavior of RWRS as $d$ increases from 4 to 5 is due to a fundamental difference between the behavior of simple random walk in 4 and 5 dimensions, a difference that is less well known than the fundamental recurrence/transience dichotomy between 2 and 3 dimensions. If we take two independent simple random walks and look at the intersection of their two trajectories, then this intersection almost surely is infinite for $1 \leq d \leq 4$ but finite for $d \geq 5$. This result is described in detail in Lawler [19], Section 3. The general necessary and sufficient condition for weak Bernoulli is that

$$|S[0, \infty) \cap S(-\infty, 0]| < \infty \text{ almost surely,}$$

i.e., the future and the past trajectories of the walk have a finite intersection, providing us with another nice dichotomy. For simple random walk, the latter condition is equivalent to $d \geq 5$, as we just described, and explains Theorem 2.5(i). This part in fact extends to a random walk with zero mean and finite variance. If, on the other hand, the random walk has a drift, then the latter condition holds because in positive time the walk is moving in the opposite direction of where it is moving in negative time. This explains Theorem 2.5(ii).

### 2.2. *Finitary factor*

We now move on to discussing when RWRS is a "finitary factor of an i.i.d. process". The following serves as a crude definition. First, if we have an i.i.d. process $W = (W_n)_{n \in \mathbb{Z}}$, taking values in a finite set $A$, and a map $f$ from $A^{\mathbb{Z}}$ to $B^{\mathbb{Z}}$, with $B$ another finite set, such that $f$ is translation invariant, then we say that the stationary random process $f(W)$ is a *factor* of $W$. If, in addition, any fixed coordinate of the image process can be determined by knowing a sufficiently large but finite (in general random) number of coordinates of the domain process relative to the position of that coordinate, then the factor map is called *finitary*.

It is well known and quite elementary to show that if a process is a finitary factor of an i.i.d. process and the so-called expected coding length (which, loosely speaking, is the expected number of bits in the domain process one needs to look at to determine a single bit of the image process) is finite, then the process is weak Bernoulli. Without the finite expected coding length assumption, this is not true. An example can be found for example in Burton and Steif [4].

As to whether being weak Bernoulli implies being a finitary factor of an i.i.d. process, del Junco and Rahe [15] constructed a counterexample. Interestingly, RWRS provides a natural counterexample. For the sake of exposition we restrict to simple random walk, but the result holds in much greater generality.

**Theorem 2.6 (Steif [39]).** *For simple random walk on $\mathbb{Z}^d$, the associated RWRS is not a finitary factor of an i.i.d. process for any $d \geq 1$.*

Observe that for $d = 1$ this result follows from Theorem 2.2. However, the proof of Theorem 2.6 is *much* simpler than the proof of Theorem 2.2.

Theorems 2.5 and 2.6 together tell us that RWRS associated with simple random walk on $\mathbb{Z}^d$, $d \geq 5$, provides us with a natural example of a random process that is weak Bernoulli but not a finitary factor of an i.i.d. process. Two key facts are used in the proof of Theorem 2.6: (i) a finitary factor of an i.i.d. process must satisfy "standard large deviation behavior", meaning that the mean ergodic theorem holds at an exponential rate (see Marton and Shields [28]); (ii) for simple random walk on $\mathbb{Z}^d$, $d \geq 1$, the number of sites visited satisfies "nonstandard large deviation behavior" (see Donsker and Varadhan [5]), and hence so does the scenery record.

It turns out that adding a drift to the random walk (which we saw causes the associated RWRS to become weak Bernoulli) results in it being a finitary factor of an i.i.d. process.

**Theorem 2.7 (Keane and Steif [18]).** *For nearest-neighbor random walk on $\mathbb{Z}$ with nonzero mean, the associated RWRS is a finitary factor of an i.i.d. process.*

We point out that this result rests on deep results in Rudolph [34], [35]. In a more direct approach, without the use of these latter papers it is possible to prove the weaker result that RWRS associated with nearest-neighbor random walk on $\mathbb{Z}$ with nonzero mean is a finitary factor of a countable state Markov chain that has "exponentially decaying return probabilities". The proof of Theorem 2.7 exploits a certain type of "regenerative structure" that is present in a random walk with positive drift. This regenerative structure is seen when observing the random walk between the time it visits $z$ for the last time and the time it visits $z + 1$ for the last time.

It is indicated in [18] how to generalize Theorem 2.7 to a nearest-neighbor random walk on $\mathbb{Z}^d$ with nonzero drift.

We summarize the results presented so far. For simple random walk on $\mathbb{Z}^d$, the associated RWRS, while K is not Bernoulli for $d = 1$ and 2. When we move to $d = 3$ and 4, it becomes Bernoulli, but not yet weak Bernoulli. When we move to $d \geq 5$, it becomes weak Bernoulli, but not yet a finitary factor of an i.i.d. process. Finally, when a drift is added, it becomes a finitary factor of an i.i.d. process.

### 2.3. Reduced RWRS

Which of the results in Sections 2.1 and 2.2 survive when we look at the second coordinate of RWRS alone, i.e., the scenery record without the steps of the walk. It is evident that all "positive results" survive the reduction. However, it is not clear which "negative results" do.

The following are two negative results that do survive the reduction, generalizing Theorems 2.2 and 2.6.

**Theorem 2.8 (Hoffman [10]).** *For simple random walk on $\mathbb{Z}$, the associated reduced RWRS is not Bernoulli.*

(Hoffman actually proves that the process is not even "loosely Bernoulli" as Kalikow did.)

**Theorem 2.9 (Steif [39]).** *For simple random walk on $\mathbb{Z}^d$, the associated reduced RWRS is not a finitary factor of an i.i.d. process for any $d \geq 1$.*

The proof of Theorem 2.8 consists of the following ingredients: (1) Kalikow's result that RWRS is not (loosely) Bernoulli; (2) Matzinger's result that "scenery reconstruction" is possible; (3) Thouvenot's "relative isomorphism theory"; (4) Rudolph's result that a "two-point weakly mixing extension" of a Bernoulli process is a Bernoulli process.

Scenery reconstruction will be described in Section 3. Ingredient (2) is crucial because it makes up for the loss of the first coordinate when going from RWRS to reduced RWRS. Indeed, for recurrent random walk the combination of the walk *and* the scenery record in RWRS allows us to retrieve the full scenery from a single realization of RWRS.

### *2.4. Non-i.i.d. random scenery or random walk*

If the i.i.d. random scenery is replaced by a stationary random field, then the situation becomes more subtle. In Section 2.1, we already mentioned that Meilijson proved Theorem 2.1 under the much weaker assumption that the random scenery is totally ergodic. In den Hollander [11], it was pointed out that Theorem 2.1 holds if and only if the random scenery is ergodic w.r.t. the subgroup of translations generated by $F - F = \{z - z' : z, z' \in F\}$. Beyond this, results are so far limited.

One early result is the following. Suppose that the random scenery is obtained by taking an irrational rotation of the circle and putting a $+1$ each time the top of the circle is hit and a $-1$ when the bottom half is hit. Suppose that the random walk stands still with probability $\frac{1}{2}$ and moves one unit to the right with probability $\frac{1}{2}$. Then, as was shown by Adler and Shields [1], [2] using geometric arguments, the associated RWRS is Bernoulli. In Shields [37] a combinatorial proof was given, and it was proved that the associated RWRS is not weak Bernoulli.

The following is a generalization of Theorem 2.5.

**Theorem 2.10 (den Hollander, Keane, Serafin and Steif [12]).**
(i) *If the random scenery is non-atomic, then for simple random walk on $\mathbb{Z}^d$, $1 \leq d \leq 4$, the associated RWRS is not weak Bernoulli.*
(ii) *If the random scenery is the plus state of the low temperature Ising model, then for simple random walk on $\mathbb{Z}^d$, $d \geq 5$, the associated RWRS is weak Bernoulli.*

To obtain Theorem 2.10(ii), one needs to be able to control the dependencies in the random field on pairs of infinite sets that are far away from each other. In the case of the low temperature Ising model, the relevant methods were developed in Burton and Steif [3] with the help of techniques from percolation. These methods can be carried over to a more general class of Markov random fields, leading to extensions of Theorem 2.10(ii) (see [12]). Van der Wal [40] has further extended Theorem 2.10(ii) to a class of random fields that are "sufficiently rapidly mixing". This class includes the two-dimensional Ising model at arbitrary supercritical temperatures, as well as all $d$-dimensional Gibbs measures at sufficiently high temperatures. He also generalized Theorem 2.3 to random sceneries that are "exponentially mixing", which includes these same sets of examples.

**Theorem 2.11 (van der Wal [40]).** *Under the same conditions as in Theorem 2.3, if the random scenery is exponentially mixing, then the associated RWRS is not Bernoulli.*

Kalikow [16] states that, once Theorem 2.2 is obtained, one can argue (using abstract ergodic theory) that the same result holds for a random scenery with positive entropy.

As far as replacing the steps of the random walk by a stationary random process is concerned, results are again limited. We mention one key result, which generalizes Theorem 2.3.

**Theorem 2.12 (Steif [39]).** *Let $X = (X_n)_{n \in \mathbb{Z}}$ be a stationary random process taking values in $\mathbb{Z}^d$ such that $S = (S_n)_{n \in \mathbb{Z}}$ is transient. If $X$ is Bernoulli and the random scenery is i.i.d., then the associated RWRS is Bernoulli.*

The technique for proving Theorem 2.12 is very different from that of proving Theorem 2.4. The former is obtained by constructing an explicit factor map from an i.i.d. process to the RWRS, while the latter is obtained by verifying the coupling property that characterizes Bernoulli for RWRS: the so-called "very weak Bernoulli" property (see [13]).

Much remains to be done in further relaxing the mixing conditions on scenery and walk.

### 2.5. Induced RWRS

Consider reduced RWRS. Suppose that we condition on the scenery at the origin being $+1$, consider only those times at which a $+1$ appears in the scenery record, and report the scenery seen by the walk at such times, i.e.,

$$\widehat{C} = (\widehat{C}_k)_{k \in \mathbb{Z}} \quad \text{with} \quad \widehat{C}_k = (C_{z + S_{T_k}})_{z \in \mathbb{Z}^d},$$

where

$$T_0 = 0,$$
$$T_k = \inf\{n > T_{k-1} \colon Y_n = +1\}, \, k \in \mathbb{N},$$
$$T_k = \sup\{n < T_{k+1} \colon Y_n = +1\}, \, k \in -\mathbb{N}.$$

The process $\widehat{C}$ is called the INDUCED RWRS (with $\{Y_0 = +1\}$ the induction set).

Clearly, $\widehat{C}$ is stationary. By Kakutani's random ergodic theorem, if $C$ is ergodic, then so is $\widehat{C}$. Mixing properties, however, are in general not inherited under induction. The following is a positive result. Here, the $\sigma$-field at infinity consists of those events that do not depend on $C_z$ for $z$ in any finite subset of $\mathbb{Z}^d$.

**Theorem 2.13 (den Hollander [11]).** *Suppose that $C$ has a trivial $\sigma$-field at infinity. If $Y_1$ is not constant almost surely, then induced RWRS is strongly mixing.*

The strong mixing property was first proved by Keane and den Hollander [17] for the case where $C$ is i.i.d. and the random walk is transient. Their proof uses a specific coupling technique, which was extended in [11] to cover the general case stated in the theorem. The coupling is delicate especially for recurrent random walk. Later Georgii [9] used a stronger form of coupling, called *orbit coupling*, weakened the condition on $C$ to it being ergodic w.r.t. the subgroup of translations generated by $F - F = \{z - z' \colon z, z' \in F\}$, and proved that under this weaker condition $\widehat{C}$ is a K-automorphism.

### 2.6. Conditional probabilities

Let us return to i.i.d. scenery and walk. We next investigate continuity properties of conditional probabilities for RWRS.

Given a general stationary random process $W = (W_n)_{n \in \mathbb{Z}}$ taking values in $\{-1, +1\}^{\mathbb{Z}}$, we may ask whether there is a version $V(\cdot \mid \eta)$ of the conditional probability distribution ($P$ is the law of $W$)

$$P\left(W_0 \in \cdot \mid W = \eta \text{ on } \mathbb{Z} \backslash \{0\}\right), \quad \eta \in \{-1, +1\}^{\mathbb{Z} \backslash \{0\}},$$

such that the map $\eta \mapsto V(\cdot \mid \eta)$ is *continuous*. If there is not, then we may ask whether there is a version that is *continuous almost everywhere*. These types of questions are of interest in probability theory and in statistical physics. Indeed, it turns out that many natural transformations acting on the set of stationary random sequences are capable of turning random sequences for which the first statement holds into random sequences for which not even the second statement holds. The history and recent developments of this research area are highlighted in the proceedings of a workshop held at EURANDOM in December 2003 (see van Enter, Le Ny and Redig [7]).

To tackle the question about the existence of "nice" conditional probabilities, a key concept is the following. For $n \in \mathbb{N}$, let $\Lambda_n = [-n, n] \cap \mathbb{Z}$. A configuration $(W_n)_{n \neq 0}$ is said to be a BAD CONFIGURATION for $W_0$ if there is an $\epsilon > 0$ such that for all $n \in \mathbb{N}$ there are $m \in \mathbb{N}$ with $m \geq n$ and $\delta \in \{-1, +1\}^{\mathbb{Z} \backslash \{0\}}$ with $\delta = \eta$ on $\Lambda_n \backslash \{0\}$ such that

$$\begin{aligned} \big\| P\left(W_0 \in \cdot \mid W = \eta \text{ on } \Lambda_m \backslash \{0\}\right) \\ - P\left(W_0 \in \cdot \mid W = \delta \text{ on } \Lambda_m \backslash \{0\}\right) \big\| \geq \epsilon, \end{aligned}$$

where $\| \cdot \|$ denotes the total variation norm. In words, by tampering with the configuration outside any given finite box, the conditional distribution of the coordinate at the origin can be nontrivially affected. A configuration that is not bad is called good.

The importance of these notions is described in Maes, Redig and Van Moffaert [27], where it is shown that every version of the conditional probability distribution is discontinuous at all the bad configurations, while there exists a version that is continuous at all the good configurations. More details can be found in den Hollander, Steif and van der Wal [14].

Returning to RWRS, the following results show that interesting behavior occurs for the conditional probability distribution of $Y_0$ (the scenery value at time 0) given $(Z_n)_{n \neq 0}$. Once more we restrict to special classes of random walks, though the results hold in much greater generality.

**Theorem 2.14 (den Hollander, Steif and van der Wal [14]).**
(i) *For arbitrary random walk, the set of bad configurations $(Z_n)_{n \neq 0}$ for $Y_0$ is non-empty.*

(ii) *For simple random walk on $\mathbb{Z}^d$, there is a version of the conditional probability distribution of $Y_0$ given $(Z_n)_{n\neq 0}$ that is continuous almost everywhere if and only if $d = 1$ or 2.*

If instead we consider the conditional probability distribution of $X_0$ (the step at time 0) given $(Z_n)_{n\neq 0}$, then the answer is different.

**Theorem 2.15 (den Hollander, Steif and van der Wal [14]).**
(i) *For arbitrary random walk, the set of bad configurations $(Z_n)_{n\neq 0}$ for $X_0$ is non-empty.*
(ii) *For simple random walk on $\mathbb{Z}^d$, there is a version of the conditional probability distribution of $X_0$ given $(Z_n)_{n\neq 0}$ that is continuous almost everywhere if and only if $1 \leq d \leq 4$.*

Theorems 2.14 and 2.15 rely on a *full classification* of the bad configurations. What drives these results are the same intersection properties of the random walk that are behind Theorem 2.5. We point out that, in Theorems 2.14 and 2.15, for the good configurations the random variable at the origin is determined by a sufficiently large piece of the good configuration. For example, any realization where the walker eventually returns to the origin allows us to read off the scenery value at the origin. Since the latter cannot change with any further information on the process after the return, this explains Theorem 2.14(ii).

If we consider the conditional probability distribution of $Y_0$ given $(Y_n)_{n\neq 0}$ instead, i.e., we pose the continuity problem for reduced RWRS, then we believe that very different behavior occurs.

**Conjecture 2.16.** *Consider the random walk on $\mathbb{Z}$ with $F = \{-1, +1\}$, $m_F(+1) = p$ and $m_F(-1) = 1 - p$, where $p \in [\frac{1}{2}, 1]$. For $p \in [\frac{1}{2}, \frac{4}{5})$ every configuration $(Y_n)_{n\neq 0}$ is bad for $Y_0$, while for $p \in (\frac{4}{5}, 1]$ every configuration $(Y_n)_{n\neq 0}$ is good for $Y_0$.*

We are presently attempting to prove this conjecture. If true, it would provide us with a remarkable example where the continuity problem has a *phase transition* in the drift parameter.

## 3. Scenery reconstruction

Scenery reconstruction is the problem of recovering $C$ given only $Y_+ = (Y_n)_{n\in\mathbb{N}_0}$. In words, given a *single* forward realization of the scenery record (i.e., the scenery as seen through the eyes of the walker at nonnegative times), is it possible to reconstruct the *full* scenery *without* knowing the walk? Remarkably, the answer is "sometimes yes". In this section we mention a number of results that have been obtained in past years. A detailed overview, including a description of the main techniques, is given in Lember and Matzinger [22], which also contains a full bibliography. The scenery reconstruction problem was raised by den Hollander and Keane and by Benjamini and Kesten in the mid 1980's. Most of the progress on this difficult problem has been achieved only since the late 1990's.

A precise formulation of the scenery reconstruction problem requires that we make some fair restrictions:

1. Scenery reconstruction is not possible when the random walk is transient, because then almost surely there are sites that the walker never visits.
2. For simple random walk, any scenery that puts +1's in the interval $\Lambda_k$, $k \in \mathbb{N}$, cannot be distinguished from the scenery that is obtained by shifting it $l$ units to the left or to the right with $l \leq k$ even. Similarly, if the random walk is reflection invariant, then any two sceneries that are reflections of each other cannot be distinguished.

3. Lindenstrauss [23] has constructed a countably infinite collection of one-dimensional sceneries (all different under translation and reflection) that cannot be distinguished with simple random walk. The sceneries in this collection have a certain "self-similar structure" and therefore have measure zero under i.i.d. scenery processes.

Thus, scenery reconstruction is at best possible for *recurrent* random walk, *up to translation and reflection* (in general), for *almost every* scenery and *almost surely* w.r.t. the walk. It turns out that these restrictions are enough, and so the problem can now be formalized, as follows.

Recall that $F$ and $G$ are the sets of possible values of the walk increments, respectively, the scenery components. Two sceneries $C$ and $C'$ are said to be *equivalent*, written $C \sim C'$, if they can be obtained from each other by a translation or reflection. Then scenery reconstruction is said to be possible when there exists a measurable function $A\colon G^{\mathbb{N}_0} \to G^{\mathbb{Z}^d}$, called a RECONSTRUCTION ALGORITHM, such that

$$P(A(Y_+) \sim C) = 1.$$

Several methods have been developed to deal with scenery reconstruction in different cases. These methods vary substantially with modifications of $F$ and $G$. Roughly speaking, the larger $F$ is, the harder the scenery reconstruction, while the larger $G$ is, the easier the scenery reconstruction. See Lember and Matzinger [22] for insight into why.

Most results so far are restricted to i.i.d. random walk and i.i.d. random scenery. We mention three key results, all for $d = 1$ and for $m_F, m_G$ the uniform distribution on $F, G$.

**Theorem 3.1 (Matzinger [29], [30]).** *Scenery reconstruction is possible when $F = \{-1, +1\}$ and $|G| = 3$, or $F = \{-1, 0, +1\}$ and $|G| = 2$.*

**Theorem 3.2 (Löwe, Matzinger and Merkl [26]).** *Scenery reconstruction is possible when $|F| < |G|$.*

**Theorem 3.3 (Lember and Matzinger [20], [21]).** *Scenery reconstruction is possible when $F \supsetneq \{-1, 0, +1\}$ and $|G| = 2$.*

(The uniformity restrictions on $m_F$ and $m_G$ may be relaxed.)

Matzinger and Rolles [32] have shown that a finite piece of scenery can be reconstructed in a time that is polynomial in the length of the piece. In addition, Matzinger and Rolles [31] have shown that scenery reconstruction is robust against errors: if the scenery record $Y_+$ is perturbed by randomly changing each bit with a probability $\epsilon > 0$, then the scenery can still be reconstructed from the perturbed scenery record, provided $\epsilon$ is small enough.

Scenery reconstruction is particularly challenging in two dimensions. This is because recurrent random walk on $\mathbb{Z}^2$ returns to sites extremely slowly. Very little is known so far, the most notable result being due to Löwe and Matzinger [24], stating that scenery reconstruction is possible for simple random walk on $\mathbb{Z}^2$ when $|G|$ is sufficiently large. The proof requires $|G|$ to be very large.

Very little progress has been made so far on scenery reconstruction for non-i.i.d. sceneries. One result can be found in Löwe and Matzinger [25]. Partial progress is underway for one-dimensional Gibbs sceneries by Lember and Matzinger (private communication).

## 4. Open questions

We close by formulating a number of open questions:

1. Does the technical condition in Theorem 2.3 hold for arbitrary random walk, as conjectured in [13]?
2. Is there an analogue of Theorem 2.4 for transient random walk and non-i.i.d. random scenery similar in spirit to Theorem 2.11, which generalizes Theorem 2.3 to recurrent random walk and exponentially mixing random sceneries?
3. How far can Theorem 2.10(ii) be extended within the class of Gibbs random sceneries, in particular, in the non-uniqueness regime?
4. Let $X = (X_n)_{n \in \mathbb{Z}}$ be a stationary random process taking values in $\mathbb{Z}$ with zero mean, implying that the random walk $S = (S_n)_{n \in \mathbb{Z}}$ is recurrent (Durrett [6], Section 6.3). For i.i.d. scenery, under what conditions on $X$ is the associated RWRS not Bernoulli? This would generalize both Theorem 2.2 and Theorem 2.12? Can it be Bernoulli?
5. When is the induced RWRS Bernoulli or weak Bernoulli? This would generalize Theorem 2.13 and its subsequent extension in [9].
6. Is scenery reconstruction possible as soon as the entropy of $m_F$ is strictly smaller than that of $m_G$, as conjectured in [26]? This would generalize Theorem 3.2.
7. Is scenery reconstruction possible for an arbitrary recurrent random walk with i.i.d. increments and an arbitrary i.i.d. random scenery?
8. To what extent can scenery reconstruction be carried through for non-i.i.d. random sceneries, e.g. Gibbs random sceneries?

A final reference on RWRS of interest is Gantert, König and Shi [8], where the small, moderate and large deviation behavior of sums of scenery values seen along the walk are reviewed. This paper includes many references to the relevant literature.